\documentclass[11pt]{amsart}

\usepackage{amsmath}
\usepackage{amssymb}
\usepackage{mathrsfs}
\usepackage[all]{xy}

\newtheorem{theorem}{Theorem}[section]
\newtheorem{claim}[theorem]{Claim}

\newtheorem{lemma}[theorem]{Lemma}

\newtheorem{corollary}[theorem]{Corollary}
\newtheorem{obs}[theorem]{Observation}

\theoremstyle{definition}
\newtheorem{definition}[theorem]{Definition}

\newtheorem{question}[theorem]{Question}

\theoremstyle{remark}

\newcount\skewfactor
\def\mathunderaccent#1#2 {\let\theaccent#1\skewfactor#2
\mathpalette\putaccentunder}
\def\putaccentunder#1#2{\oalign{$#1#2$\crcr\hidewidth
\vbox to.2ex{\hbox{$#1\skew\skewfactor\theaccent{}$}\vss}\hidewidth}}
\def\name{\mathunderaccent\tilde-3 }

\def\smallbox#1{\leavevmode\thinspace\hbox{\vrule\vtop{\vbox
   {\hrule\kern1pt\hbox{\vphantom{\tt/}\thinspace{\tt#1}\thinspace}}
   \kern1pt\hrule}\vrule}\thinspace}

\newcommand{\cf}{{\rm cf}}

\def\qedref#1{$\qed_{\reforiginal{#1}}$}

\title{Cardinal characteristics at aleph omega}
\author{Shimon Garti}
\address{Einstein Institute of Mathematics,
 The Hebrew University of Jerusalem,
 Jerusalem 91904, Israel}
\email{shimon.garty@mail.huji.ac.il}

\author{Moti Gitik}
\address{School of Mathematical Sciences,
 Raymond and Beverly Sackler Faculty of Exact Science,
 Tel Aviv University,
 Ramat Aviv 69978, Israel}
\email{gitik@post.tau.ac.il}

\author{Saharon Shelah}
\address{Institute of Mathematics
 The Hebrew University of Jerusalem,
 Jerusalem 91904, Israel
 and  Department of Mathematics
 Rutgers University
 New Brunswick, NJ 08854, USA}
\email{shelah@math.huji.ac.il}
\urladdr{http://www.math.rutgers.edu/\char`\~shelah}
\thanks{The research of the first and third authors was supported by European Research Council, Grant no. 338821. Second author's research was partially supported by Israel Science Fountdation Grant no. 1216/18. This is publication 1156 of the third author}


\begin{document}
\let\labeloriginal\label
\let\reforiginal\ref
\def\ref#1{\reforiginal{#1}}
\def\label#1{\labeloriginal{#1}}

\begin{abstract}
We prove the consistency of $\mathfrak{u}_{\aleph_\omega}<2^{\aleph_\omega}$.
We also show that the consistency strength of this statement is the existence of a measurable cardinal $\kappa$ with $o(\kappa)=\kappa^{++}$.
\end{abstract}

\maketitle

\newpage

\section{Introduction}

The ultrafilter number is one of the cardinal characteristics defined on the continuum. 
The definition generalizes to uncountable cardinals, and makes sense over both regular and singular cardinals.
Reacll that a filter $\mathscr{F}$ over $\lambda$ is called \emph{uniform} iff every element of $\mathscr{F}$ is of size $\lambda$.

\begin{definition}
\label{defuf} The ultrafilter number. \newline 
Let $\lambda$ be an infinite cardinal, and $\mathcal{F}$ a uniform filter over $\lambda$.
\begin{enumerate}
\item [$(\aleph)$] A base $\mathcal{A}$ for $\mathcal{F}$ is a subfamily of $\mathcal{F}$ such that for every $B\in\mathcal{F}$ there is some $A\in\mathcal{A}$ with the property $A\subseteq^*B$.
\item [$(\beth)$] The character of $\mathcal{F}$ is the minimal size of a base for $\mathcal{F}$, denoted by ${\rm Ch}(\mathcal{F})$.
\item [$(\gimel)$] The ultrafilter number $\mathfrak{u}_\lambda$ is the minimal cardinality of a filter base for some uniform ultrafilter over $\lambda$.
\end{enumerate}
\end{definition}

An easy diagonalization argument shows that $\mathfrak{u}_\lambda>\lambda$ for every infinite cardinal $\lambda$, see \cite[Claim 1.2]{MR2992547}.
Typically, $\mathfrak{u}_\lambda$ tends to be large.
Our main objective in this paper is to prove the consistency of $\mathfrak{u}_\lambda<2^\lambda$ where $\lambda=\aleph_\omega$.
We shall obtain, further, the inequality $\mathfrak{u}_\lambda<\mathfrak{d}_\lambda$ at $\aleph_\omega$, and show that this implies $\mathfrak{u}_\lambda<\mathfrak{i}_\lambda$ as well.
Another important issue is the consistency strength of these results.
We shall prove that a measurable cardinal $\kappa$ such that $o(\kappa)=\kappa^{++}$ is the exact consistency strength.

The consistency of $\mathfrak{u}_\lambda<2^\lambda$ for a strong limit singular cardinal $\lambda$ has been proved in \cite{MR2992547}, starting from a supercompact cardinal in the ground model.
Let us try to explain the difficulty of obtaining this result at $\aleph_\omega$.
Suppose that $\lambda>\cf(\lambda)=\kappa$ is a limit of measurable cardinlas.
The basic idea in \cite{MR2992547} is to choose a sequence of measurable cardinals $\langle\lambda_i:i\in\kappa\rangle$ such that $\lambda=\bigcup_{i\in\kappa}\lambda_i$ and $2^{\lambda_i}=\lambda_i^+$ for every $i\in\kappa$.
We fix a normal ultrafilter $\mathscr{U}_i$ over $\lambda_i$ for every $i\in\kappa$, and using the local instances of GCH at the $\lambda_i$s we arragne a base of $\mathscr{U}_i$ of size $\lambda_i^+$ whose elements form a $\subseteq^*$-decreasing sequence.
This requires the normality of the ultrafilters.

By forcing an appropriate true cofinality for the products of $\lambda_i$s and $\lambda_i^+$s, one can define a uniform ultrafilter $\mathscr{U}$ over $\lambda$, which is basically the sum of the $\mathscr{U}_i$s with respect to some uniform ultrafilter over $\kappa$.
One can collapse this configuration to $\aleph_\omega$, while keeping the pcf structure.
However, normal ultrafilters do not exist over the $\aleph_n$s.
Since normality plays a key-role in the above construction, it is not clear whether the inequality $\mathfrak{u}_\lambda<2^\lambda$ is obtainable where $\lambda=\aleph_\omega$, and this has been asked in \cite{MR3832086}.

The main idea in the second section of this paper is that one can replace ultrafilters by filters with some strong properties.
We can keep the normality of our filters, and this can be done without the expensive requirement of being an ultrafilter.
In particular, the objects needed for proving $\mathfrak{u}_\lambda<2^\lambda$ may live over small accessible cardinals, thus giving the consistency of $\mathfrak{u}_{\aleph_\omega}<2^{\aleph_\omega}$.
Further results about cardinal characteristics are proved at strong limit singular carinals, including the case of $\aleph_\omega$.

In the third section we try to analyze the consistency strength of the above result.
In order to force the failure of SCH, as done for our result, one needs at least a measurable cardinal $\kappa$ such that $o(\kappa)=\kappa^{++}$ as proved in \cite{MR1098782}.
We shall see that this is the exact consistency strength of $\mathfrak{u}_{\aleph_\omega}<2^{\aleph_\omega}$.
Finally, we pose in the last section some open problems concerning characteristics at strong limit singular cardinals.

Our notation is mostly standard.
For background in pcf theory and Prikry-type forcing notions we suggest \cite{MR2768693} and \cite{MR2768695} respectively.
We shall use the Jerusalem forcing notation, so if $p,q$ are conditions in $\mathbb{P}$ then $p\leq q$ means that $q$ is stronger than $p$.
We denote the usual L\'evy collapse by ${\rm Col}(\kappa,\lambda)$ or ${\rm Col}(\kappa,<\lambda)$.

\newpage

\section{The ultrafilter number and other animals}

A quasi order $(W,\leq_W)$ is a transitive and reflexive binary relation.
A subset $V$ of $W$ is dense iff for every $x\in W$ there is some $y\in V$ such that $x\leq_W y$.
The set $V$ is open iff it is $\leq_W$-upward closed.
The definition of \emph{sullam} below is a generalization of the concept of scale from pcf theory.

\begin{definition}
\label{defsullam} Sullam. \newline 
Let $\kappa$ be a regular cardinal, $J$ an ideal over $\kappa$ such that $J^{\rm bd}_\kappa\subseteq J, \langle W_i:i\in\kappa\rangle$ a sequence of quasi orders and $\bar{f}=\langle f_\alpha:\alpha\in\lambda\rangle$ a sequence of functions where $f_\alpha\in\prod_{i\in\kappa}W_i$ for each $\alpha\in\lambda$. \newline 
We shall say that $\bar{f}$ is a sullam in $(\prod_{i\in\kappa}W_i,J)$ iff the following requirements are met:
\begin{enumerate}
\item [$(\aleph)$] $\bar{f}$ is $J$-increasing, that is $\alpha<\beta\Rightarrow \{i\in\kappa:f_\alpha(i)\leq_{W_i}f_\beta(i)\}=\kappa\ {\rm mod}\ J$.
\item [$(\beth)$] $\bar{f}$ is cofinal in the sense that if $V_i$ is a dense subset of $W_i$ for every $i\in\kappa$ then there exists an ordinal $\alpha\in\lambda$ such that $\{i\in\kappa:f_\alpha(i)\in V_i\}=\kappa\ {\rm mod}\ J$.
\end{enumerate}
\end{definition}

Let $\mathscr{D}$ be a filter over a regular cardinal $\theta$.
Throughout the paper we will always assume that our filters are uniform.
Moreover, we always assume that such $\mathscr{D}$ extends the collection $\{\theta-u:u\in[\theta]^{<\theta}\}$.
If $\mathscr{D}$ is a filter over $\theta$ then $I(\mathscr{D}) = \{A\subseteq\theta:\theta-A\in\mathscr{D}\}$ is the dual ideal.
The collection of $\mathscr{D}$-positive sets is $\mathscr{D}^+ = \mathcal{P}(\theta)-I(\mathscr{D})$.
If $A,B\subseteq\theta$ then $A\subseteq_{\mathscr{D}}B$ iff $A-B\in I(\mathscr{D})$.

Let $\mathscr{D}$ be a filter over $\theta, W$ a quasi order and $g:W\rightarrow\mathscr{D}^+$ a function.
We shall say that $g$ is $\subseteq_{\mathscr{D}}$-decreasing iff $s\leq_W t \Rightarrow g(t)\subseteq_{\mathscr{D}}g(s)$.
We shall say that $g$ has the decidability property iff for every $s\in W$ and every $A\subseteq\theta$ there exists some $t\in W$ such that $s\leq_W t$ and $(g(t)\subseteq_{\mathscr{D}}A)\vee(g(t)\subseteq_{\mathscr{D}}(\theta-A))$.
The idea is to apply this property to small accessible cardinals.
On these cardinals there are no noraml ultrafilters.
However, there are normal filters, and if we force decidability over these filters then we can imitate the effect of an ultrafilter.

\begin{definition}
\label{defnice} Nice systems. \newline 
Assume that $\kappa=\cf(\mu)<\mu$. \newline 
A nice system $\mathcal{S}$ for $\mu$ consists of the following objects:
\begin{enumerate}
\item [$(a)$] An increasing sequence of cardinals $(\mu_i:i\in\kappa)$ such that $\mu=\bigcup_{i\in\kappa}\mu_i$.
\item [$(b)$] An increasing sequence of regular cardinals $(\lambda_i:i\in\kappa)$ such that $\mu_i\leq\lambda_i<\mu_{i+1}$ for every $i\in\kappa$.
\item [$(c)$] A sequence $(\mathscr{D}_i:i\in\kappa)$ of filters, each $\mathscr{D}_i$ is defined over $\lambda_i$.
\item [$(d)$] A sequence $(W_i:i\in\kappa)$ of quasi orders.
\item [$(e)$] A sequence of functions $(g_i:i\in\kappa)$ such that $g_i:W_i\rightarrow\mathscr{D}_i^+$ is $\subseteq_{\mathscr{D}_i}$-decreasing and every $g_i$ has decidability.
\end{enumerate}
\end{definition}

The main theorem of this section produces ultrafilters with small bases over a singular cardinal from nice systems over this cardinal.
In the next section we shall discuss how to force the existence of nice systems.

\begin{theorem}
\label{thmsmallu} Let $\mu>\cf(\mu)=\kappa$ be a singular cardinal and let $\lambda=\cf(\lambda)\in(\mu,2^\mu]$. \newline 
Let $\mathcal{S}$ be a nice system for $\mu$, and assume that:
\begin{enumerate}
\item [$(\aleph)$] $\mathscr{D}$ is a uniform ultrafilter over $\kappa$, generated by $\leq\lambda$ many sets.
\item [$(\beth)$] $\bar{f}=(f_\alpha:\alpha\in\lambda)$ is a sullam in $(\prod_{i\in\kappa}W_i,\mathscr{D})$.
\item [$(\gimel)$] $\cf(\prod_{i\in\kappa}\mathscr{D}_i,\supseteq)=\lambda$.
\end{enumerate}
Then there exists a uniform ultrafilter $\mathscr{U}$ over $\mu$ such that ${\rm Ch}(\mathscr{U})\leq\lambda$.
\end{theorem}

\par\noindent\emph{Proof}. \newline 
We define our ultrafilter over $\mu$ as follows.
Given $A\subseteq\mu$ we shall say that $A\in\mathscr{U}$ iff there exist an ordinal $\alpha\in\lambda$ and an element $x\in\mathscr{D}$ such that $i\in x \Rightarrow g_i(f_\alpha(i))\subseteq_{\mathscr{D}_i}(A\cap\lambda_i)$.
We claim that $\mathscr{U}$ is a uniform ultrafilter over $\mu$ and ${\rm Ch}(\mathscr{U})\leq\lambda$.
It follows from the definition that $\mathscr{U}$ is a filter, and from the decidability properties of the functions $g_i$s we infer that $\mathscr{U}$ is an ultrafilter.
Notice that $\mathscr{U}$ is uniform since if $A\in\mathscr{U}$ then the inclusion $g_i(f_\alpha(i))\subseteq_{\mathscr{D}_i}(A\cap\lambda_i)$ means that $A\cap\lambda_i$ is of size $\lambda_i$ and this happens at $\kappa$ many $i$s as $\mathscr{D}$ is uniform, and hence $|A|=\mu$.
We must show, therefore, that ${\rm Ch}(\mathscr{U})\leq\lambda$.

To see this, fix a sequence $(\mathcal{A}_\alpha:\alpha\in\lambda)$ of elements of $\prod_{i\in\kappa}\mathscr{D}_i$ which is cofinal.
Explicitly, $\mathcal{A}_\alpha = \langle A_{\alpha i}:i\in\kappa\rangle$ for every $\alpha\in\lambda$, $A_{\alpha i}\in\mathscr{D}_i$ and if $\langle B_i:i\in\kappa\rangle\in \prod_{i\in\kappa}\mathscr{D}_i$ then for some $\alpha\in\lambda$ we have $\{i\in\kappa:A_{\alpha i}\subseteq B_i\} =_\mathscr{D} \kappa$.
The existence of $(\mathcal{A}_\alpha:\alpha\in\lambda)$ is ensured by $(\gimel)$.

Enumerate the elements of $\mathscr{D}$ by $(x_\gamma:\gamma\in\lambda)$ (repetitions are welcome).
For every triple $(\alpha,\beta,\gamma)$ such that $\alpha,\beta,\gamma\in\lambda$ let $B_{\alpha\beta\gamma} = \bigcup\{g_i(f_\alpha(i))\cap A_{\beta i}:i\in x_\gamma\}$.
Define:
$$
\mathcal{B} = \{B_{\alpha\beta\gamma}:\alpha,\beta,\gamma\in\lambda\}.
$$
We claim that $\mathscr{U}$ is generated by $\mathcal{B}$.
First of all, notice that each $B_{\alpha\beta\gamma}$ belongs to $\mathscr{U}$ by the choice of $(\mathcal{A}_\alpha:\alpha\in\lambda)$ and the definition of $\mathscr{U}$.
Second, assume that $A\in\mathscr{U}$.
From the definition of $\mathscr{U}$ we can fix an element $x\in\mathscr{D}$ and an ordinal $\alpha\in\lambda$ such that $g_i(f_\alpha(i))\subseteq_{\mathscr{D}_i}(A\cap\lambda_i)$ for every $i\in x$.
Let $\gamma\in\lambda$ be such that $x=x_\gamma$.

For each $i\in x_\gamma$ let $B_i\in\mathscr{D}_i$ be an element which satisfies $g_i(f_\alpha(i))\cap B_i\subseteq A\cap\lambda_i$.
For each $i\in\kappa-x_\gamma$ let $B_i$ be $\lambda_i$.
The sequence $\langle B_i:i\in\kappa\rangle$ is an element of $\prod_{i\in\kappa}\mathscr{D}_i$.
By virtue of $(\gimel)$ we can find an ordinal $\beta\in\lambda$ such that $i\in\kappa\Rightarrow A_{\beta i}\subseteq B_i$.
It follows that $g_i(f_\alpha(i))\cap A_{\beta i}\subseteq g_i(f_\alpha(i))\cap B_i\subseteq A\cap\lambda_i$ for every $i\in x_\gamma$.
This means that $B_{\alpha\beta\gamma}\subseteq\bigcup_{i\in x_\gamma}(A\cap\lambda_i)=A$, so we are done.

\hfill \qedref{thmsmallu}

Remark that if $\mu$ is a strong limit cardinal which is our main interest, then $2^\kappa<\mu<\lambda$ and hence every ultrafilter $\mathscr{D}$ over $\kappa$ is generated by less than $\lambda$ many sets.
Observe also that actually ${\rm Ch}(\mathscr{U})=\lambda$ in the above construction, a fact which can be proved as done in \cite{MR3832086}.
By the forcing construction of the next section it follows that one can increase $2^{\aleph_\omega}$ to some regular $\tau<\aleph_{\omega_4}$ and obtain ${\rm Sp}_\chi(\aleph_\omega)\supseteq{\rm Reg}\cap[\aleph_{\omega+1},\tau]$.
This gives a positive answer to a question from \cite{MR3832086}.

It has been proved in \cite{dear} that $\mathfrak{u}_\lambda<\mathfrak{d}_\lambda$ is consistent for some strong limit singular cardinal $\lambda$, where $\mathfrak{d}_\lambda$ is the minimal size of a dominating family of functions from $\lambda$ into $\lambda$.
The forcing construction of the next section shows that this can be forced over $\lambda=\aleph_\omega$ as well.
Based on this observation, we wish to prove the consistency of $\mathfrak{u}_\lambda<\mathfrak{i}_\lambda$ for $\lambda=\aleph_\omega$, by proving that $\mathfrak{d}_\lambda\leq\mathfrak{i}_\lambda$ whenever $\lambda>\cf(\lambda)=\omega$.

\begin{definition}
\label{defind} Independence. \newline 
Assume that $\lambda>\cf(\lambda)=\omega$.
\begin{enumerate}
\item [$(\aleph)$] For a family $\mathcal{A} = \{I_\alpha:\alpha\in\kappa\}\subseteq[\lambda]^\lambda$ let ${\rm comb}(\mathcal{A})$ be the collection of all finite boolean combinations of elements of $\mathcal{A}$. Explicitly, an element of ${\rm comb}(\mathcal{A})$ is a set of the form $\bigcap\Gamma - \bigcup\Delta$ where $\Gamma,\Delta\in[\mathcal{A}]^{<\omega}$ and $\Gamma\cap\Delta=\varnothing$.
\item [$(\beth)$] A family $\mathcal{A}\subseteq[\lambda]^\lambda$ is independent iff ${\rm comb}(\mathcal{A})\subseteq[\lambda]^\lambda$.
\item [$(\gimel)$] The independence number $\mathfrak{i}_\lambda$ is the minimal size of a maximal independent family in $[\lambda]^\lambda$.
\end{enumerate}
\end{definition}

The cardinal characteristic $\mathfrak{i}_\lambda$ is a generalization of the independence number $\mathfrak{i}$.
One has to be careful here, since the generalization to $\mathfrak{i}_\kappa$ for uncountable $\kappa$ in the sense of boolean combinations of size less than $\kappa$ is problematic.
In particular, it is not clear whether such a maximal family exists.
It is possible to speak about finite boolean combinations, but this seems less natural.
However, in the case of $\lambda>\cf(\lambda)=\omega$, where the natural requirement is boolean combinations of size less than $\cf(\lambda)$, we are speaking again about finite combinations, and the concept is well-defined.

We shall prove below that $\mathfrak{d}_\lambda\leq\mathfrak{i}_\lambda$ whenever $\lambda>\cf(\lambda)=\omega$.
Recall that $\mathfrak{d}_\lambda>\lambda$ for every infinite cardinal $\lambda$.
It follows from this inequality that $\mathfrak{i}_\lambda>\lambda$,
but we include a direct short proof to this fact.

\begin{obs}
\label{obsabovelam} If $\lambda>\cf(\lambda)=\omega$ then $\mathfrak{i}_\lambda>\lambda$.
\end{obs}

\par\noindent\emph{Proof}. \newline 
Let $\mathcal{A}=\{I_\alpha:\alpha\in\lambda\}$ be an independent family of size $\lambda$.
We shall prove that $\mathcal{A}$ is not maximal.
Call $B$ a boolean combination of $\mathcal{A}$ iff there are disjoint finite sets $\Gamma,\Delta\subseteq\lambda$ such that $B = \bigcap_{\alpha\in\Gamma}I_\alpha - \bigcup_{\alpha\in\Delta}(\lambda-I_\alpha)$.
Each boolean combination of $\mathcal{A}$ is of size $\lambda$, as $\mathcal{A}$ is independent.
We may assume without loss of generality that $\mathcal{A}$ contains all of its boolean combinations.
We may assume, further, that each element of $\mathcal{A}$ appears $\lambda$-many times in the enumeration $\{I_\alpha:\alpha\in\lambda\}$.

By induction on $\alpha\in\lambda$ choose two distinct elements $\gamma_\alpha,\delta_\alpha\in I_\alpha$ so that the following hold:
\begin{enumerate}
\item [$(a)$] $\gamma_\alpha\notin\{\gamma_\beta,\delta_\beta:\beta<\alpha\}$.
\item [$(b)$] $\delta_\alpha\notin\{\gamma_\beta,\delta_\beta:\beta<\alpha\}$.
\end{enumerate}
Define $I = \{\gamma_\alpha:\alpha\in\lambda\}$ and notice that $\{\delta_\alpha:\alpha\in\lambda\}\subseteq(\lambda-I)$.
Moreover, if $B\in{\rm comb}(\mathcal{A})$ then $|B\cap I|=|B\cap(\lambda-I)|=\lambda$ since $B$ appears $\lambda$-many times in the enumeration and by the construction.
If $\alpha\in\lambda$ then $\delta_\alpha\in I_\alpha-I$ and $\gamma_\alpha\in I_\alpha-(\lambda-I)$, so $I\notin\mathcal{A}$.
We conclude that $\mathcal{A}\cup\{I\}$ is independent and $\mathcal{A}\subsetneq\mathcal{A}\cup\{I\}$, so $\mathcal{A}$ is not maximal as required.

\hfill \qedref{obsabovelam}

Recall that an unsplittable family $\mathcal{A}\subseteq[\lambda]^\lambda$ is a collection of elements in $[\lambda]^\lambda$ such that no single element of $[\lambda]^\lambda$ splits them all.
The reaping number $\mathfrak{r}_\lambda$ is the minimal size of an unsplittable family in $[\lambda]^\lambda$.
Our goal is to prove that $\mathfrak{r}_\lambda\leq\mathfrak{i}_\lambda$ and $\mathfrak{d}_\lambda\leq\mathfrak{i}_\lambda$. 
The first inequality is easy, but we spell out the argument:

\begin{claim}
\label{clmri} Let $\lambda>\cf(\lambda)=\omega$. \newline 
Then $\mathfrak{r}_\lambda\leq\mathfrak{i}_\lambda$.
\end{claim}

\par\noindent\emph{Proof}. \newline 
Suppose that $\mathfrak{i}_\lambda=\kappa$ and $\mathcal{A} = \{I_\alpha:\alpha\in\kappa\}$ is maximal independent.
We assume that ${\rm comb}(\mathcal{A})\subseteq\mathcal{A}$.
Let us show that $\mathcal{A}$ is unsplittable.
If not, then there is some $A\in[\lambda]^\lambda$ which splits any element of $\mathcal{A}$.
This means that if $B\in\mathcal{A}$ then $|B\cap A|=|B\cap(\lambda-A)|=\lambda$, so $\mathcal{A}\cup\{A\}$ is independent.
Notice that $A\notin\mathcal{A}$ since $B\cap(\lambda-A)\neq\varnothing$ for every $B\in\mathcal{A}$, so $\mathcal{A}\cup\{A\}$ contradicts the maximality of $\mathcal{A}$.
We conclude, therefore, that $\mathcal{A}$ is unspittable, so $\mathfrak{r}_\lambda\leq\kappa=\mathfrak{i}_\lambda$ and we are done.

\hfill \qedref{clmri}

Remark that this simple argument does not imply $\mathfrak{u}_\lambda\leq\mathfrak{i}_\lambda$, see Question \ref{qudi}.
For proving that $\mathfrak{d}_\lambda\leq\mathfrak{i}_\lambda$ we need a lemma which ensures the existence of a special kind of pseudointersections.
For $A,B\in[\lambda]^\lambda$ we say that $B\subseteq^*_{\rm end}A$ iff $B$ is expressible as $\bigcup_{n\in\omega}E_n, m<n<\omega\Rightarrow|E_m|\leq|E_n|$ and there exists some $n_0\in\omega$ such that $\bigcup_{n_0\leq n\in\omega}E_n\subseteq A$.

\begin{lemma}
\label{lempseudo} Assume that:
\begin{enumerate}
\item [$(\aleph)$] $\lambda>\cf(\lambda)=\omega$.
\item [$(\beth)$] $(C_n:n\in\omega)\subseteq[\lambda]^\lambda$ is $\subseteq^*$-decreasing.
\item [$(\gimel)$] $\mathcal{A}\subseteq[\lambda]^\lambda, |\mathcal{A}|<\mathfrak{d}_\lambda$.
\item [$(\daleth)$] For every $n\in\omega$ and every $A\in\mathcal{A}$ it is true that $|A\cap C_n|=\lambda$.
\end{enumerate}
Then there exists a set $B\in[\lambda]^\lambda$ such that:
\begin{enumerate}
\item [$(a)$] $B=\bigcup_{n\in\omega}E_n$, and $m<n<\omega\Rightarrow|E_m|\leq|E_n|$.
\item [$(b)$] For every $m\in\omega, B\subseteq^*_{\rm end}C_m$.
\item [$(c)$] For every $A\in\mathcal{A}$ and every $m\in\omega, |\bigcup_{m\leq n\in\omega}E_n\cap A|=\lambda$.
\end{enumerate}
\end{lemma}

\par\noindent\emph{Proof}. \newline 
By cutting each $C_n$ with the previous elements of the sequence we may assume that $(C_n:n\in\omega)$ is $\subseteq$-decreasing.
Fix an increasing sequence of ordinals $(\alpha_n:n\in\omega)$ such that $\lambda = \bigcup_{n\in\omega}\alpha_n$.
For every $A\in\mathcal{A}$ and every $m\in\omega$ define $f_A\in{}^\lambda\lambda$ by letting $f_A(\beta)$ be the $\beta$-th element of the set $C_{n(\beta)}\cap A$, where $n(\beta)\in\omega$ is the unique natural number such that $\alpha_n\leq\beta<\alpha_{n+1}$.

Since $|\mathcal{A}|<\mathfrak{d}_\lambda$, we can choose a function $g\in{}^\lambda\lambda$ such that $\forall A\in\mathcal{A}, \neg(g\leq^* f_A)$.
This means that the set $\{\beta\in\lambda:f_A(\beta)<g(\beta)\}$ is of size $\lambda$, and it happens at every $A\in\mathcal{A}$.
For each $n\in\omega$, define $E_n = E_n^g = \bigcup\{C_{n(\beta)}\cap g(\beta): \alpha_n\leq\beta<\alpha_{n+1}\}$.
Set $B = \bigcup_{n\in\omega}E_n^g$.
We claim that $B$ exemplifies the lemma.
Notice that $(a)$ follows directly from the construction, so we are left with the other two statements.

For $(b)$ fix $m\in\omega$.
If $n>m$ then any element $\gamma\in E_n^g$ satisfies $\gamma\in C_m$ since the sequence $(C_n:n\in\omega)$ is $\subseteq$-decreasing, and due to the definition of $E_n^g$.
It follows that $E_n\subseteq C_m$ for every $n>m$, so $(b)$ is satisfied.
For $(c)$ fix $A\in\mathcal{A}, m\in\omega$.
The set $\{\beta\in\lambda:f_A(\beta)<g(\beta)\}$ is of cardinality $\lambda$, hence also the set $\{\beta\in(\alpha_{m+1},\lambda):f_A(\beta)<g(\beta)\}$ is of size $\lambda$.
If $\beta$ belongs to this set then $f_A(\beta)\in E_{n(\beta)}\cap A$ and since $f_A$ is one-to-one we infer that $|\bigcup_{n\geq m}E_n\cap A|=|(\alpha_{m+1},\lambda)|=\lambda$.

\hfill \qedref{lempseudo}

We shall use the above lemma within the proof of the following:

\begin{theorem}
\label{thmdi} Assume that $\lambda>\cf(\lambda)=\omega$. \newline 
Then $\mathfrak{d}_\lambda\leq\mathfrak{i}_\lambda$.
\end{theorem}

\par\noindent\emph{Proof}. \newline 
Let $\mathcal{I}$ be an independent family such that $|\mathcal{I}|<\mathfrak{d}_\lambda$.
We shall try to find a set in $[\lambda]^\lambda$ which can be added to $\mathcal{I}$ while keeping its independence, thus proving the theorem.
As a first stage we choose a countable set $\{D_n:n\in\omega\}\subseteq\mathcal{I}$ and we let $\mathcal{J} = \mathcal{I}-\{D_n:n\in\omega\}$.
Let $\mathcal{A}$ be ${\rm comb}(\mathcal{J})$, and remark that $|\mathcal{A}|<\mathfrak{d}_\lambda$.

For every $y\in[\lambda]^\lambda$ let $y^0=y, y^1=\lambda-y$.
Let $\mathscr{X}$ be the topological space ${}^\omega 2$ with the product topology, where each component is given the discrete topology.
For each $f\in\mathscr{X}$ and every $n\in\omega$ let $C_n = C_n^f = \bigcap_{m<n}D_m^{f(m)}$.
Observe that for each $f\in\mathscr{X}$, the sequence $(C_n^f:n\in\omega)$ is $\subseteq$-decreasing.
Observe also that each $C_n^f$ intersects each $A\in\mathcal{A}$ at $\lambda$ many points, so we can apply Lemma \ref{lempseudo} to get $B_f\in[\lambda]^\lambda$ such that:
\begin{enumerate}
\item [$(a)$] $B_f = \bigcup_{n\in\omega}E_n^f$.
\item [$(b)$] For every $n\in\omega, B_f\subseteq^*_{\rm end}C_n^f$.
\item [$(c)$] For every $A\in\mathcal{A}, m\in\omega$ we have $|\bigcup_{n\geq m}E_n^f\cap A|=\lambda$.
\end{enumerate}
We make the comment that if $f,g\in\mathscr{X}$ and $f\neq g$ then there is some $m = m(f,g)\in\omega$ such that $\bigcup_{n\geq m}E_n^f\cap \bigcup_{n\geq m}E_n^g = \varnothing$.
Indeed, if $m$ is the first point in which $f(m)\neq g(m)$ then $D_m^{f(m)}\cap D_m^{g(m)}=\varnothing$.
Now we use the fact that $(C_n:n\in\omega)$ is $\subseteq$-decreasing and the definition of the $E_n$s in the lemma.

Since $\mathscr{X}$ is separable we can choose two disjoint countable dense sets $P,Q$. Enumerate the elements of $P\cup Q$ by $\{f_n:n\in\omega\}$.
By induction on $n\in\omega$ we omit from $B_{f_n}$ the following set.
If $f_n\in P$ then for any $\ell<n$ such that $f_\ell\in Q$ we let $m_\ell = m(f_\ell,f_n)$ and we omit from $B_{f_n}$ the set $\bigcup_{k\leq m_\ell}E_k^{f_\ell}$.
Similarly, if $f_n\in Q$ then for any $\ell<n$ such that $f_\ell\in P$ we let $m_\ell = m(f_\ell,f_n)$ and we omit from $B_{f_n}$ the set $\bigcup_{k\leq m_\ell}E_k^{f_\ell}$.
This process keeps $(a),(b),(c)$ above, and results in a family $\{B_{f_n}:n\in\omega\}$ whose elements are pairwise disjoint.
Define $Y = \bigcup_{f\in P}B_f$ and $Z = \bigcup_{f\in Q}B_f$, and notice that $Y\cap Z=\varnothing$.

We claim that $\mathcal{I}\cup\{Y\}$ is an independent family.
In order to prove this, fix a couple of finite sets $s,t\subseteq\mathcal{I}$.
Denote $u = s\cap\mathcal{J}, v = t\cap\mathcal{J}$.
Choose a sufficiently large $n\in\omega$ such that if $D_k\in s\cup t$ then $k<n$. This is possible since $s,t$ are finite.
Choose $f\in P$ such that if $D_k\in s$ then $f(k)=0$ and if $D_k\in t$ then $f(k)=1$.
The choice is possible since there are only finitely many $D_k$s in $s\cup t$, and setting the values of $f$ on this finite collection defines a basic open set in $\mathscr{X}$.
Since $P$ is dense in $\mathscr{X}$ we can pick up any $f\in P$ which belongs to this open set.
Consider the following statements:
\begin{enumerate}
\item [$(\alpha)$] $\bigcap s - \bigcup t = (\bigcap u - \bigcup v)\cap \bigcap\{D_k^{f(k)}: D_k\in s\cup t\}$. \newline 
This statement follows simply from the definition of these sets.
\item [$(\beta)$] $(\bigcap u - \bigcup v)\cap \bigcap\{D_k^{f(k)}: D_k\in s\cup t\} \supseteq (\bigcap u - \bigcup v)\cap \bigcap\{D_k^{f(k)}: k<n\}$. \newline 
This statement follows from the choice of $n$.
\item [$(\gamma)$] $|B_f\cap(\bigcap s - \bigcup t)|=\lambda$. \newline 
For this fact, choose $m\in\omega$ such that $\bigcup_{n\geq m}E_n^f\subseteq \bigcap_{k<n}D_k^{f(k)} = C_n^f$, using property $(b)$ of $B_f$.
From property $(c)$ we have $|\bigcup_{n\geq m}E_n^f\cap (\bigcap u - \bigcup v)|=\lambda$ and hence $|B_f\cap C_n^f\cap (\bigcap u - \bigcup v)|=\lambda$.
\item [$(\delta)$] $|Y\cap(\bigcap s - \bigcup t)|=\lambda$. \newline 
For this fact, recall that $B_f\subseteq Y$.
\end{enumerate}
By an identical argument we can show that $|Z\cap(\bigcap s - \bigcup t)|=\lambda$, upon replacing $P$ by $Q$.
This means that $\mathcal{I}\cup\{Y\}$ is independent, and the proof is accomplished.

\hfill \qedref{thmdi}

Combining the above theorem with our ability to decrease the ultrafilter number and concomitantly increase the dominating number, we obtain the following conclusion:

\begin{corollary}
\label{corui} Assuming sufficiently large cardinals in the ground model one can force $\mathfrak{u}_{\aleph_\omega}<\mathfrak{i}_{\aleph_\omega}$.
\end{corollary}

\hfill \qedref{corui}

\newpage

\section{Consistency strength}

In this section we try to determine the consistency strength of the statement $\mathfrak{u}_\lambda<2^\lambda$ where $\lambda$ is a strong limit singular cardinal.
Since $\mathfrak{u}_\lambda\geq\lambda^+$, the statement $\mathfrak{u}_\lambda<2^\lambda$ implies in particular $\lambda^+<2^\lambda$ so SCH fails at $\lambda$.
Hence it follows from \cite{MR1098782} that the consistency strength of $\mathfrak{u}_\lambda<2^\lambda$ at a strong limit singular cardinal $\lambda$ is at least a measurable cardinal $\mu$ such that $o(\mu)=\mu^{++}$.

The main result of this section is that actually this is the exact consistency strength.
For simplicity, we shall show how to obtain $\mathfrak{u}_\lambda<2^\lambda$ from $o(\mu)=\mu^{++}+1$, and then explain how to modify the proof in order to get the desired consistency strength of $o(\mu)=\mu^{++}$.

The forcing notion to be used in this section is a relative of the extender-based Prikry forcing with interleaved collapses.
We shall define it using the notation of \cite{MR2768695}.
Let $\mu$ be a measurable cardinal and let $E$ be a $(\mu,\lambda)$-extender over $\mu$.
Let $\jmath_E:V\rightarrow M\cong{\rm Ult}(V,E)$ be the canonical embedding.
We are assuming GCH in the ground model.

\begin{definition}
\label{defbroken} Extender-based Prikry with broken collapses. \newline 
Let $A$ be an unbounded subset of measurable cardinals below $\mu$, and let $h:\mu\rightarrow\mu$ be the increasing enumeration of the elements of $A$.
We define the extender-based Prikry forcing $\mathbb{P}$ with interleaved collapses broken at $A$. \newline 
A condition $p\in\mathbb{P}$ has the form:
\begin{center}
$\{\langle 0,\langle\tau_1,\ldots,\tau_n\rangle,\langle f_0,\ldots,f_n\rangle,F\} \cup$ \newline  $\{\langle\gamma,p^\gamma\rangle:\gamma\in{\rm supp}(p)-\{0,{\rm mc}(p)\}\} \cup$ \newline $\{\langle{\rm mc}(p),p^{\rm mc(p)},T\rangle\}$.
\end{center}
The part which consists of $\{\langle 0,\langle\tau_1,\ldots,\tau_n\rangle\} \cup \{\langle\gamma,p^\gamma\rangle:\gamma\in{\rm supp}(p)-\{0,{\rm mc}(p)\}\} \cup \{\langle{\rm mc}(p),p^{\rm mc(p)},T\rangle\}$ is the usual condition of the extender-based Prikry forcing.
Each $p^\gamma$ is an approximation to the Prikry sequence at the $\gamma$th coordinate, ${\rm mc}(p)$ is the maximal coordinate of the condition $p$, and $T$ is the tree associated with ${\rm mc}(p)$. \newline 
The sequence $\langle f_0,\ldots,f_n\rangle$ and the function $F$ describe the collapses, and should satisfy the following:
\begin{itemize}
\item $f_0\in{\rm Col}(\omega,\tau_1)$.
\item $f_i\in{\rm Col}(\tau_i^{+3},h(\tau_i))\times {\rm Col}(h(\tau_i)^+,\tau_{i+1})$ for every $i\in(0,n)$.
\item $f_n\in{\rm Col}(\tau_n^{+3},h(\tau_n))\times {\rm Col}(h(\tau_n)^+,\mu)$.
\item $F$ is a function defined on the projection of $T$ to the normal ultrafilter $E_\mu$ and satisfies $F(\langle\nu_0,\ldots,\nu_{i-1}\rangle)\in{\rm Col}(\nu_{i-1}^{+3},<h(\nu_{i-1}))\times {\rm Col}(h(\nu_{i-1})^+,\mu)$.
\end{itemize}
\end{definition}

Ahead of depicting the orders $\leq_{\mathbb{P}}$ and $\leq_{\mathbb{P}}^*$ we indicate that the only deviation in our version from the traditional extender-based Prikry forcing with interleaved collapses is the definition of the collapses.
In the traditional version each $f_i$ belongs to ${\rm Col}(\tau_i^{+3},\tau_{i+1})$ while in the current version we decompose the collapse and require $f_i\in{\rm Col}(\tau_i^{+3},h(\tau_i))\times {\rm Col}(h(\tau_i)^+,\tau_{i+1})$.
The reason is that we wish to preserve \emph{measurability-like properties} of the $h(\tau_i)$s.
Of course, full measurability will not be preserved since these cardinals are going to be below $\aleph_\omega$ in the generic extension.
However, appropriate filters with decidability will be obtained.

\begin{definition}
\label{deforders} The orders $\leq_{\mathbb{P}}$ and $\leq_{\mathbb{P}}^*$. \newline 
Assume that $p,q\in\mathbb{P}$ where $\mathbb{P}$ is the extender-based Prikry forcing with interleaved collapses broken at $A$.
Let $n_p,n_q$ denote $\ell g(p^0),\ell g(q^0)$ respectively. \newline 
We shall say that $p\leq_{\mathbb{P}}q$ iff:
\begin{enumerate}
\item [$(a)$] $\{\langle 0,q^0\rangle\}\cup \{\langle\gamma,q^\gamma\rangle: \gamma\in{\rm supp}(q) - \{0,{\rm mc}(q)\}\}\cup \{{\rm mc}(q),q^{\rm mc},T^q\}$ extends $\{\langle 0,p^0\rangle\}\cup \{\langle\gamma,p^\gamma\rangle: \gamma\in{\rm supp}(p) - \{0,{\rm mc}(p)\}\}\cup \{{\rm mc}(p),p^{\rm mc},T^p\}$ by the usual order of the extender-based Prikry forcing.
\item [$(b)$] For every $i<n_p, f_i^p\leq f_i^q$.
\item [$(c)$] For every $\eta\in T^{q,0}_{q^0}, F^p(\eta)\subseteq F^q(\eta)$.
\item [$(d)$] For every $i\in[n_p,n_q), F^p((q^0-p^0)\upharpoonright(i+1))\subseteq f^q_i$.
\item [$(e)$] $\sup({\rm rang}(f_{n_p}))<\min(q^0-p^0)$.
\end{enumerate}
For $p,q\in\mathbb{P}$ we shall say that $p\leq_{\mathbb{P}}^* q$ iff $p\leq_{\mathbb{P}}q$ and for every $\gamma\in{\rm supp}(p), p^\gamma=q^\gamma$.
\end{definition}

The basic properties of this forcing notion are similar to the properties of the traditional extender-based Prikry forcing with interleaved collapses.
In particular, $\mathbb{P}$ satisfies the $\mu^{++}$-cc and has the Prikry property. The cardinal $\mu^+$ remains a cardinal in the generic extension as well. For detailed proofs see \cite{MR2768695}.
Let $G\subseteq\mathbb{P}$ be a generic set over $V$.
In $V[G]$ we have $2^{\aleph_n}=\aleph_{n+1}$ for every $n\in\omega$ and $2^{\aleph_\omega}=\aleph_{\omega+2}$.

\begin{theorem}
\label{thmo+1} Let $\mu$ be a measurable cardinal such that $o(\mu)=\mu^{++}+1$. \newline 
Then there is a generic extension which satisfies the following:
\begin{enumerate}
\item [$(\aleph)$] $\mu=\aleph_\omega$.
\item [$(\beth)$] $2^{\aleph_\omega}=\aleph_{\omega+2}$.
\item [$(\gimel)$] If $\theta<\aleph_\omega$ or $\theta\geq\aleph_\omega^+$ then $2^\theta=\theta^+$.
\item [$(\daleth)$] There exists a nice system $\mathcal{S}$ for $\aleph_\omega$ which satisfies the assumption of Theorem \ref{thmsmallu} with $\kappa=\omega$ and $\lambda=\aleph_{\omega+1}$.
\end{enumerate}
\end{theorem}

\par\noindent\emph{Proof}. \newline 
We begin with GCH in the ground model, and we fix a $(\mu,\mu^+)$-extender $E$.
Let $\jmath:V\rightarrow M\cong{\rm Ult}(V,E)$ be the associated ultrapower embedding, so $\mu={\rm crit}(\jmath)$ and $M\supseteq V_{\mu+2}$.
Recall that the first non-trivial measure $E_\mu = \{A\subseteq\mu:\mu\in\jmath(A)\}$ is a noraml ultrafilter over $\mu$.
Let $\jmath_\mu:V\rightarrow M_\mu\cong{\rm Ult}(V,E_\mu)$ be the ultrapower embedding by $E_\mu$.
Let $k_\mu:M_\mu\rightarrow M$ be the induced elementary embedding, namely $k_\mu([f]_{E_\mu})=\jmath(f)(\mu)$.
It is easy to check that the following diagram commutes:

$$
\xymatrix{
{V} \ar[ddrr]_{\jmath_\mu} \ar[rr]^\jmath & & M \\ \\ 
& & M_\mu \ar[uu]_{k_\mu} }
$$

Since $E_\mu\in\mathcal{P}(\mathcal{P}(\mu))$ and $M\supseteq V_{\mu+2}$ we see that $E_\mu\in M$.
It follows that $\mu$ is a measurable cardinal in $M$.
Moreover, the set $\{\nu\in\mu:\nu$ is a measurable cardinal$\}$ is unbounded in $\mu$.
By elementarity we see that in $M_\mu$ the set $\{\nu\in\jmath_\mu(\mu):\nu$ is a measurable cardinal$\}$ is unbounded in $\jmath_\mu(\mu)$.

Since $(\mu^{++})^{M_\mu}<\jmath_\mu(\mu)$ we can fix in $M_\mu$ a measurable cardinal $\eta$ such that $M_\mu\models\mu^{++}<\eta<\jmath_\mu(\mu)$.
Let $\mathscr{D}_\eta$ be a normal ultrafilter over $\eta$ in $M_\mu$.
Elements of $M_\mu$ are represented by functions from $\mu$ into $V$, so pick a function $h:\mu\rightarrow\mu$ such that $\jmath_\mu(h)(\mu)=\eta$ and another function $r:\mu\rightarrow V_\mu$ such that $\jmath_\mu(r)(\mu)=\mathscr{D}_\eta$.

Let $\mathbb{P}$ be the extender-based Prikry forcing with interleaved collapses broken at $A$ where $A={\rm rang}(h)$.
We define, in $V[G]$, our nice system $\mathcal{S}$.
Let $(\mu_i:i\in\omega)$ be the Prikry sequence associated with the normal measure $E_\mu$.
For every $i\in\omega$ let $\lambda_i = h(\mu_i)$, and $\mathscr{D}_i = r(\mu_i)$.
Without loss of generality, each $\lambda_i$ is measurable in $V$ and $\mathscr{D}_i$ is a normal ultrafilter over $\lambda_i$.

For each $i\in\omega$, the part of the forcing above $\lambda_i$ does not change the measurability of $\lambda_i$.
Similarly, the collapse below $\mu_i$ does not affect the measurability of $\lambda_i$.
The part of ${\rm Col}(\mu_i^{+3},<\lambda_i)$ will destroy the measurability of $\lambda_i$.
Nevertheless, $\mathscr{D}_i$ will generate in $V[G]$ a precipitous normal filter over $\lambda_i$, and by abuse of notation let us call this filter $\mathscr{D}_i$.
For every $i\in\omega$ let $\jmath_i:V\rightarrow{\rm Ult}(V,r(\mu_i))$ be the canonical embedding and let $W_i = {\rm Col}(\mu_i^{+3},<\jmath_i(h(\mu_i)))$.
Notice that each $W_i$ is $\mu_i$-complete as a partial order.
Likewise, the partial order $(\mathscr{D}_i^+,\supseteq)$ is isomorphic, as a forcing notion, to $W_i$ for every $i\in\omega$.
Fix an isomorphism $g_i:W_i\rightarrow\mathscr{D}_i^+$ for every $i\in\omega$.
This accomplishes the definition of the nice system, and all the requirements are easily verified.

Let $\mathscr{D}$ be a uniform ultrafilter over $\omega$.
We shall prove in the lemmata below that there exists a sullam of length $\lambda=\aleph_{\omega+1}$ in $(\prod_{i\in\omega}W_i,\mathscr{D})$ and that $\cf(\prod_{i\in\kappa}\mathscr{D}_i,\supseteq)=\lambda$.
This proves part $(\daleth)$ of the theorem.
The other parts follow from the properties of the extender-based Prikry forcing with interleaved collapses, so we are done.

\hfill \qedref{thmo+1}

Recall that $\bar{f}$ is a sullam in $(\prod_{i\in\omega}W_i,\mathscr{D})$ iff $\bar{f}$ is $\mathscr{D}$-increasing and $\mathscr{D}$-cofinal with respect to the dense subsets of the quasi orders $W_i$.

\begin{lemma}
\label{lemsullam} In the generic extension defined in the above theorem there exists a sullam $\bar{f}=\langle f_\alpha:\alpha\in\aleph_{\omega+1}\rangle$ in $(\prod_{i\in\omega}W_i,\mathscr{D})$.
\end{lemma}

\par\noindent\emph{Proof}. \newline 
Let $\mu=\aleph_\omega$.
The main point in the proof is a statement which says that one can construct a sequence $\langle\bar{V}_\xi:\xi\in\mu^+\rangle$ with the following properties:
\begin{enumerate}
\item [$(a)$] $\bar{V}_\xi=\langle V_{\xi i}:i\in\omega\rangle\in \prod_{i\in\omega}W_i$ for every $\xi\in\mu^+$.
\item [$(b)$] $V_{\xi i}$ is open and dense in $W_i$ for every $i\in\omega, \xi\in\mu^+$.
\item [$(c)$] If $U_i$ is a dense open subset of $W_i$ for every $i\in\omega$ then there exists an ordinal $\xi\in\mu^+$ and a natural number $i_0\in\omega$ such that $V_{\xi i}\subseteq U_i$ for every $i\in[i_0,\omega)$.
\end{enumerate}
Before proving this statement let us show how to derive the conclusion of the lemma from the existence of such a sequence.

We define $f_\alpha\in\prod_{i\in\omega}W_i$ by induction on $\alpha\in\mu^+$.
If $\alpha=0$ then we choose $w^0_i\in V_{0i}\subseteq W_i$ for every $i\in\omega$ and simply define $f_0(i)=w^0_i$ for every $i\in\omega$.
If $\alpha=\beta+1$ and $f_\beta$ is at hand then we choose $w^\alpha_i\in V_{\alpha i}\subseteq W_i$ such that $w^\beta_i\leq_{W_i}w^\alpha_i$ for every $i\in\omega$.
This can be done since $V_{\alpha i}$ is dense.
Finally, let $\alpha\in\mu^+$ be a limit ordinal, so $\cf(\alpha)<\mu$.
Let $(f_{\alpha_\gamma}:\gamma\in\cf(\alpha))$ be a sequence of functions such that $(\alpha_\gamma:\gamma\in\cf(\alpha))$ is cofinal in $\alpha$.
Choose $i_0\in\omega$ such that $\cf(\alpha)<\mu_{i_0}$.
Recall that $W_i$ is $\mu_i$-complete for every $i\in\omega$.
Hence, if $i\in[i_0,\omega)$ then one can choose $w^\alpha_i\in V_{\alpha i}\subseteq W_i$ such that $w^{\alpha_\beta}_i\leq_{W_i}w^\alpha_i$ for every $\beta<\cf(\alpha)$.
By letting $f_\alpha(i)=w^\alpha_i$ for each $i\in[i_0,\omega)$ and zero otherwise we complete the construction of the functions, and the sequence $(f_\alpha:\alpha\in\mu^+)$ is as required.

We are left with the construction of $\langle\bar{V}_\xi:\xi\in\mu^+\rangle$.
For this end, we shall prove the following claim which is a bit stronger than the property that we need.
We claim that if $U_i$ is a dense open subset of $W_i$ for every $i\in\omega$ then there exists a function $h:V_\mu\rightarrow V_\mu$ in the ground model, such that the range of $h$ is a dense open subset of $U_i$ for every $i\in\omega$.
This claim implies the statement at the beginning of the proof by our cardinal arithmetic assumptions.
Indeed, the amount of all functions $h:V_\mu\rightarrow V_\mu$ in the ground model is $\mu^+$.
By collecting all the functions whose range is a dense subset of $W_i$ for every $i\in\omega$ and enumerating them as $\langle\bar{V}_\xi:\xi\in\mu^+\rangle$ we will be done.

For proving the claim we work in $V$, so let $(\name{U}_i:i\in\omega)$ be a name for a sequence $(U_i:i\in\omega)$ of dense open subsets of the $W_i$s, and assume that this is forced by the weakest condition.
In $V$ we have a name $\name{W}_i$ for the collapse ${\rm Col}(\name{\mu}_i^{+3},<\name{\tau}_i)$, and we shall use the fact that $\name{W}_i$ is forced to be $\name{\mu}_i^{+3}$-complete.
Suppose, therefore, that $(\mu_j:j<i)$ has been decided, and we try to consider in $V$ all the possibilities for members of $\name{W}_i$.

Let $p$ be a condition in $\mathbb{P}$.
By taking a direct extension of $p$ if needed, we may assume that elements of $\name{W}_i$ are decided independently from the stage above $i$, see \cite{MR1007865}.
Let $A_p$ be the measure one set associated with $p$ at the $i$th stage.
For every $\nu\in A_p$ let $\nu^0$ be its projection to $E_\mu$, the normal ultrafilter of the extender $E$.
Fix an element $\nu\in A_p$ and its projection $\nu^0$.
Notice that $\nu^0$ points to the collapse ${\rm Col}((\nu^0)^{+3},<\tau_{\nu^0})$.
Enumerate the elements of this collapse in some canonical way which depends only on $\nu^0$.

Let $x_0$ be the first element according to this enumeration.
Define $B_p(\nu^0) = \{\sigma\in A_p:\sigma^0=\nu^0\}$.
Upon enumerating the elements of $B_p(\nu^0)$ we choose for every $\sigma\in B_p(\nu^0)$ some $t_\sigma$ such that the following hold:
\begin{enumerate}
\item [$(a)$] There exists a condition $q$ such that $p^\frown\nu\leq^*q$.
\item [$(b)$] $q\Vdash x_0\leq_{W_i}t_\sigma\in\name{U}_i$.
\item [$(c)$] $q\Vdash t_{\sigma'}\leq_{W_i}t_\sigma$ for every $\sigma'\in B_p(\nu^0)$ which appeared before $\sigma$ in our enumeration.
\end{enumerate}
The choice is possible since $|B_p(\nu^0)|\leq(\nu^0)^{++}$ and the pertinent collapse is $(\nu^0)^{+3}$-complete.
Moreover, we can choose at the end of the inductive process an element $y_0(\nu^0)\in{\rm Col}((\nu^0)^{+3},<\tau_{\nu^0})$ such that $x_0\leq_{W_i} t_\sigma\leq_{W_i}y_0(\nu^0)$ for every $\sigma\in B_p(\nu^0)$. It follows that if $q\geq p$ and the $i$th element of the Prikry sequence through $E_\mu$ determined by $q$ is $\nu^0$ then $q\Vdash y_0(\nu^0)\in\name{U}_i$.
We render this process and find $y_\beta(\nu^0)\in{\rm Col}((\nu^0)^{+3},<\tau_{\nu^0})$ such that $x_\beta\leq_{W_i}y_\beta(\nu^0)$ for every $x_\beta$ in the collapse.

Let $S_{\nu^0} = \{y_\beta(\nu^0):\beta<\tau_{\nu^0}\}$.
Notice that every $q\geq p$ such that $\nu^0$ is the $i$th element of the Prikry sequence according to $q$ forces $S_{\nu^0}\subseteq\name{U}_i$ and $S_{\nu^0}$ is dense in $W_i$.
Set $h(\mu_0,\ldots,\mu_{i-1},\nu^0)$ as the upward closure of $S_{\nu^0}$.
By doing this for every $\nu^0$ such that $\nu\in A_p$ we define $h:V_\mu\rightarrow V_\mu$ in the ground model.
The collection of all these functions enumerated as $\langle\bar{V}_\xi:\xi\in\mu^+\rangle$ satisfies our claim, so we are done.

\hfill \qedref{lemsullam}

Finally, we compute the cofinality of $(\prod_{i\in\omega}\mathscr{D}_i,\supseteq)$.

\begin{lemma}
\label{lemcof} In the generic extension defined in the above theorem we have $\cf(\prod_{i\in\omega}\mathscr{D}_i,\supseteq)=\aleph_{\omega+1}$.
\end{lemma}

\par\noindent\emph{Proof}. \newline 
For every $i\in\omega, \mathscr{D}_i$ is a normal ultrafilter over $\lambda_i$.
Since $2^{\lambda_i}=\lambda_i^+$ in $V[G]$, we can fix a base $\mathcal{A}_i = \{A_{\alpha i}:\alpha\in\lambda_i^+\}\subseteq\mathscr{D}_i$ for every $i\in\omega$.
As each $\mathscr{D}_i$ is normal, we can replace $A_{\alpha i}$ by $B_{\alpha i}=\Delta\{A_{\beta i}:\beta\in\alpha\}\in\mathscr{D}_i$ for every $i\in\omega,\alpha\in\lambda_i^+$.
Let $\mathcal{B}_i$ be the family $\{B_{\alpha i}:\alpha\in\lambda_i^+\}$ for every $i\in\omega$.
Notice that $\mathcal{B}_i$ is a $\subseteq^*$-decreasing base of $\mathscr{D}_i$ for every $i\in\omega$.

As proved in \cite{MR2768695}, ${\rm tcf}(\prod_{i\in\omega}\lambda_i,J) = {\rm tcf}(\prod_{i\in\omega}\lambda_i^+,J) = \aleph_{\omega+1}$, where $J=J^{\rm bd}_\omega$.
Choose a scale $(f_\zeta:\zeta\in\aleph_{\omega+1})$ in the product $\prod_{i\in\omega}\lambda_i$, and another scale $(g_\xi:\xi\in\aleph_{\omega+1})$ in $\prod_{i\in\omega}\lambda_i^+$.
For every $\zeta,\xi\in\aleph_{\omega+1}$ define:
$$
\mathcal{A}_{\zeta\xi} = \langle B_{g_\xi(i)i} - f_\zeta(i):i\in\omega\rangle.
$$
Notice that each $\mathcal{A}_{\zeta\xi}$ belongs to $\prod_{i\in\omega}\mathscr{D}_i$.
We claim that $\{\mathcal{A}_{\zeta\xi}:\zeta,\xi\in\aleph_{\omega+1}\}$ is a cofinal subset of $(\prod_{i\in\omega}\mathscr{D}_i,\supseteq)$.

To see this, let $(S_i:i\in\omega)\in\prod_{i\in\omega}\mathscr{D}_i$.
For every $i\in\omega$ choose an ordinal $\alpha_i\in\lambda_i^+$ such that $B_{\alpha_i i}\subseteq^* S_i$.
The function $h(i)=\alpha_i$ defined for every $i\in\omega$ belongs to $\prod_{i\in\omega}\lambda_i^+$, so we can choose an ordinal $\xi\in\aleph_{\omega+1}$ such that $h<^*g_\xi$.
This means that for some $i_0\in\omega$ if $i\in[i_0,\omega)$ then $B_{g_\xi(i)i}\subseteq^*S_i$.
Hence for every $i\in[i_0,\omega)$ there is an ordinal $\beta_i\in\lambda_i$ such that $B_{g_\xi(i)i}-S_i\subseteq\beta_i$.

Define $f(i)=0$ whenever $i<i_0$ and $f(i)=\beta_i$ if $i\in[i_0,\omega)$.
It follows that $f\in\prod_{i\in\omega}\lambda_i$.
Choose an ordinal $\zeta\in\aleph_{\omega+1}$ such that $f<^*f_\zeta$.
Now choose some $i_1\in[i_0,\omega)$ such that $B_{g_\xi(i)i}-f_\zeta(i)\subseteq S_i$ for every $i\in[i_1,\omega)$.
By the previous definitions, $(S_i:i\in\omega)$ is $\mathscr{D}$-covered by $\mathcal{A}_{\zeta\xi}$ in the sense of $\supseteq$, since $J^{\rm bd}_\omega\subseteq\mathscr{D}$.
The proof of the lemma is accomplished.

\hfill \qedref{lemcof}

Having these lemmata at hand, we know how to obtain $\mathfrak{u}_{\aleph_\omega}<2^{\aleph_\omega}$ from a measurable cardinal $\mu$ with $o(\mu)=\mu^{++}+1$ in the ground model.
It is possible to refine this assumption.

\begin{theorem}
\label{thmo++} Suppose that $o(\mu)=\mu^{++}$. \newline 
Then one can force $\mathfrak{u}_{\aleph_\omega}=\aleph_{\omega+1} <2^{\aleph_\omega}=\aleph_{\omega+2} = \mathfrak{d}_{\aleph_\omega}$ with GCH holds everywhere apart from $\aleph_\omega$.
\end{theorem}

\par\noindent\emph{Proof}. \newline 
We shall use the idea of \cite{MR1007865} in which the failure of SCH has been forced from a measurable cardinal $\mu$ with $o(\mu)=\mu^{++}$ in the ground model.
Based on this assumption one constructs a $(\mu,\mu^{++})$-extender $E$ by changing cofinalities and producing long Rudin-Keisler increasing sequences.
The only deviation from this model in our construction is that we keep an untouched set of measurable cardinals unbounded below $\mu$.
The role of these cardinals is to carry normal precipitous filters in the generic extension after the collapses.
These filters will satisfy the assumption of Theorem \ref{thmsmallu}, thus proving the statement $\mathfrak{u}_{\aleph_\omega}=\aleph_{\omega+1}$ in the generic extension.

So fix a $(\mu,\mu^{++})$-extender $E$ and let $\jmath_E:V\rightarrow M\cong {\rm Ult}(V,E)$ be the canonical embedding.
As usual, let $E_\mu$ denote the normal ultrafilter designated by the first coordinate of $E$.
Choose an element $A\in E_\mu$ and a function $H:A\rightarrow\mu$ such that the following requirements hold for every $\nu\in A$:
\begin{enumerate}
\item [$(a)$] $\nu$ is an inaccessible cardinal.
\item [$(b)$] $H(\nu)=\nu^{++}$.
\item [$(c)$] $\jmath_E(H)(\mu)=\mu^{++}$.
\item [$(d)$] There exists a measurable cardinal $\tau_\nu$ such that $H(\nu)<\tau_\nu<\min(A-(\nu+1))$.
\item [$(e)$] $\mathscr{D}_\nu$ is a normal ultrafilter over $\tau_\nu$.
\end{enumerate}
Notice that if $\nu_0\in A-(\nu+1)$ then $H(\nu)=\nu^{++}<\nu_0$ since $\nu_0$ is inaccessible.
We force now with the extender-based Prikry forcing $\mathbb{P}$ with broken interleaved collapses, with respect to the measurable cardinals $\tau_\nu$ for every $\nu\in A$.

Fix a generic set $G\subseteq\mathbb{P}$, and let $(\mu_n:n\in\omega)$ be the Prikry sequence through $E_\mu$ as computed by $G$.
Since we use broken collapses, we force eventually with ${\rm Col}(\mu_n^{+3},<\tau_{\mu_n})\times{\rm Col}(\tau_{\mu_n}^{+3},<\mu_{n+1})$.
The purpose of breaking the collapses, as explicated above, is to maintain nice properties of the normal ultrafilters $\mathscr{D}_{\tau_{\mu_n}}$ in the generic extension.

For every $i\in\omega$ let $\lambda_i=\tau_{\mu_i}$ and let $\mathscr{D}_i$ be the filter generated by $\mathscr{D}_{\tau_{\mu_i}}$ in $V[G]$.
Notice that $\mathscr{D}_i$ is a normal precipitous filter over $\lambda_i$.
For each $i\in\omega$ let $W_i={\rm Col}(\mu_i^{+3},<\jmath_i(\tau_{\mu_i}))$ where $\jmath_i:V\rightarrow{\rm Ult}(V,\mathscr{D}_{\tau_{\mu_i}})$ is the canonical embedding.
Remark that the forcing notion $(\mathscr{D}_i^+,\supseteq)$ is isomorphic to $W_i$ for every $i\in\omega$, hence one can fix an isomorphism (of forcing notions) $g_i:W_i\rightarrow\mathscr{D}_i^+$ for every $i\in\omega$.
The objects constructed so far form a nice system $\mathcal{S}$.

Fix a uniform ultrafilter $\mathscr{D}$ over $\omega$.
Invoking the above lemmata we see that there exists a sullam of length $\aleph_{\omega+1}$ in the product $(\prod_{i\in\omega}W_i,\mathscr{D})$ and that $\cf(\prod_{i\in\omega}\mathscr{D}_i,\supseteq)=\aleph_{\omega+1}$.
We conclude, therefore, that $\mathfrak{u}_{\aleph_\omega}=\aleph_{\omega+1}$ in $V[G]$.
Notice that ${\rm tcf}(\prod_{i\in\omega}\mu_i^{++},J^{\rm bd}_\omega) = \aleph_{\omega+2} = 2^{\aleph_\omega}$ in $V[G]$, and the same is true for ${\rm tcf}(\prod_{i\in\omega}\mu_i^{+3},J^{\rm bd}_\omega)$.
From \cite{dear} we deduce that $\mathfrak{d}_{\aleph_\omega}=\aleph_{\omega+2}$, so we are done.

\hfill \qedref{thmo++}

\newpage 

\section{Open problems}

In this section we collect several natural questions which emerged from our study.
The first question is based on the impression that it is quite easy to construct a uniform ultrafilter $\mathscr{U}$ over a singular cardinal $\lambda$ such that ${\rm Ch}(\mathscr{U})=\lambda^+$.
In particular, the consistency of $\mathfrak{u}_{\aleph_\omega}=\aleph_{\omega+1}$ shows that no large cardinal assumption is required below $\lambda$ at the generic extension.

\begin{question}
\label{qplus} Is it consistent that $\lambda$ is a strong limit singular cardinal and $\mathfrak{u}_\lambda>\lambda^+$?
\end{question}

Lest $\kappa=\cf(\kappa), \mathfrak{u}_\kappa=2^\kappa$ and $2^\kappa$ is arbitrarily large obtains easily.
For example, if one adds $\tau$ many Cohen subsets of $\kappa$ where $\tau=\cf(\tau)$ to a model of $2^\kappa=\kappa^+$ then $\mathfrak{u}_\kappa=2^\kappa=\tau$ in the generic extension.
Something similar can be done at a singular cardinal $\lambda$ if we give up strong limitude.

\begin{question}
\label{qlargeu} Is it consistent that $\lambda$ is a strong limit singular cardinal, $2^\lambda$ is large and $\mathfrak{u}_\lambda=2^\lambda$?
\end{question}

Of course, there are pcf limitations on the value of $2^\lambda$ where $\lambda>\cf(\lambda)$ is strong limit, and the above question is raised under these constraints.
A more general way to phrase similar questions is by trying to give a good description of the character spectrum of the ultrafilters over some cardinal.
Recall that ${\rm Sp}_\chi(\lambda)$ is the set of all ${\rm Ch}(\mathscr{U})$ where $\mathscr{U}$ is a uniform ultrafilter over $\lambda$.
As we have seen, one can incorporate the whole interval $[\lambda^+,2^\lambda]$ in the character spectrum of a singular cardinal $\lambda$.

\begin{question}
\label{qspec} Assume that $\lambda$ is a strong limit singular cardinal. \newline 
Is it consistent that ${\rm Sp}_\chi(\lambda)$ is perforated, i.e. for some triple of regular cardinals $\kappa_0,\kappa_1,\kappa_2\in[\lambda^+,2^\lambda]$ so that $\kappa_0<\kappa_1<\kappa_2$ we have $\kappa_0,\kappa_2\in{\rm Sp}_\chi(\lambda)$ and $\kappa_1\notin{\rm Sp}_\chi(\lambda)$?
\end{question}

In order to give a positive answer one has to prove somehow an opposite statement of Theorem \ref{thmsmallu}.
Namely, one has to show that if $\kappa\in{\rm Sp}_\chi(\lambda)$ then for some sequence of regular cardinals below $\lambda$ whose $\kappa$ is the true cofinality there are products with the properties of the assumptions in Theorem \ref{thmsmallu}.
This has to be done, probably, by focusing on some core model.
As a second step one has to force elements in ${\rm Reg}\cap[\lambda^+,2^\lambda]$ for which such sequence of regular cardinals below $\lambda$ is absent.

In the constructions defined so far, ${\rm Ch}(\mathscr{U})$ is regular.
Indeed, the value of ${\rm Ch}(\mathscr{U})$ is determined by true cofinalities of products of regular cardinals.
One may wonder about the possibility of singular cardinals in ${\rm Sp}_\chi(\lambda)$.
Using independent families, one can always create a uniform ultrafilter $\mathscr{U}$ such that ${\rm Ch}(\mathscr{U})=2^\lambda$, so by forcing $2^\lambda$ to be singular we can introduce a singular cardinal into the spectrum.
Recall that if $\lambda=\aleph_0$ or even if $\lambda$ is regular and uncountable then it is even possible to force $\mathfrak{u}_\lambda$ to be a singular cardinal, by the same reasoning.

\begin{question}
\label{qsingspec} Let $\lambda$ be a strong limit singular caridnal.
\begin{enumerate}
\item [$(\aleph)$] Is it consistent that ${\rm Ch}(\mathscr{U})<2^\lambda$ and ${\rm Ch}(\mathscr{U})$ is singular for some uniform ultrafilter $\mathscr{U}$ over $\lambda$?
\item [$(\beth)$] Is it consistent that the cofinality of ${\rm Ch}(\mathscr{U})$ is less than $\lambda$ for some uniform ultrafilter $\mathscr{U}$ over $\lambda$?
\item [$(\gimel)$] Is it consistent that $\mathfrak{u}_\lambda$ is singular?
\end{enumerate}
\end{question}

Let us move to the other side of the coin as reflected in $\mathfrak{d}_\lambda$ and $\mathfrak{i}_\lambda$.
It follows from \cite{dear} and Theorem \ref{thmdi} of the current paper that $\mathfrak{d}_\lambda$, and hence $\mathfrak{i}_\lambda$, tend to be large.
Indeed, sequences with large true cofinalities increase $\mathfrak{d}_\lambda$ and unlike $\mathfrak{u}_\lambda$ we do not need any kind of measurability (or precipitous filters) from the elements of the sequence.

\begin{question}
\label{qd} Is it consistent that $\lambda$ is a strong limit singular cardinal and $\mathfrak{d}_\lambda<2^\lambda$?
\end{question}

One thing that we use when proving that $\mathfrak{d}_\lambda$ is large is local instances of GCH along the sequence of regular cardinals below $\lambda$.
It seems that in order to force a positive answer to the above question one has to violate GCH in a strong sense below $\lambda$.
A possible approach would be to force $\mathfrak{d}_\kappa<2^\kappa$ at some large cardinal $\kappa$ (e.g. supercompact) and then to singularize $\kappa$ in such a way that the inequality will be preserved.
If such an inequality is obtainable then one may wonder whether classical inequalities like $\mathfrak{d}<\mathfrak{i}$ or $\mathfrak{d}<\mathfrak{u}$ can be forced at singular cardinals. 

\begin{question}
\label{qudi} Suppose that $\lambda$ is a strong limit singular cardinal.
\begin{enumerate}
\item [$(\aleph)$] Is it consistent that $\mathfrak{d}_\lambda<\mathfrak{i}_\lambda$?
\item [$(\beth)$] Is it consistent that $\mathfrak{d}_\lambda<\mathfrak{u}_\lambda$?
\item [$(\gimel)$] Is it consistent that $\mathfrak{i}_\lambda<\mathfrak{u}_\lambda$?
\end{enumerate}
\end{question}

Of course, the first part of the question seems easier.
Finally, we are interested in the relationship between $\mathfrak{r}_\lambda$ and $\mathfrak{u}_\lambda$.
A beautiful theorem of Aubrey, \cite{MR2058185}, says that $\mathfrak{r}<\mathfrak{d}$ implies $\mathfrak{r}=\mathfrak{u}$.

\begin{question}
\label{qru} Assume that $\lambda$ is a strong limit singular cardinal.
\begin{enumerate}
\item [$(\aleph)$] Is it provable that $\mathfrak{r}_\lambda<\mathfrak{d}_\lambda$ implies $\mathfrak{r}_\lambda=\mathfrak{u}_\lambda$?
\item [$(\beth)$] Is it consistent that $\mathfrak{r}_\lambda<\mathfrak{u}_\lambda$?
\end{enumerate}
\end{question}

\newpage 

\bibliographystyle{amsplain}
\bibliography{arlist}

\end{document}